\documentclass[a4paper,reqno,10pt]{amsart}

\usepackage{amsmath}
\usepackage{amsthm}
\usepackage{amssymb}
\usepackage{amscd} 

\usepackage{bbm} 
\usepackage{mathrsfs} 

\usepackage{lmodern} 
\usepackage[T1]{fontenc} 
\usepackage{newtxtext,newtxmath}

\usepackage{tikz-cd}
\usepackage{pifont}

\theoremstyle{newplain}
\newtheorem*{Theorem*}{Theorem} 

\newcommand{\R}{\mathbb{R}} 

\newcommand{\Rn}{\R^n}

\newcommand{\calA}{\mathscr{A}}
\newcommand{\calB}{\mathscr{B}}
\newcommand{\calC}{\mathscr{C}}

\newcommand{\calE}{\mathscr{E}}
\newcommand{\calF}{\mathscr{F}}
\newcommand{\calG}{\mathscr{G}}
\newcommand{\calH}{\mathscr{H}}
\newcommand{\calI}{\mathscr{I}}

\newcommand{\calL}{\mathscr{L}}
\newcommand{\calM}{\mathscr{M}}
\newcommand{\calN}{\mathscr{N}}

\newcommand{\calP}{\mathscr{P}}

\newcommand{\sfMSNsp}{\mathsf{MSN_{sp}}}
\newcommand{\sfLOC}{\mathsf{LOC}}
\newcommand{\sfMSN}{\mathsf{MSN}}
\newcommand{\sfLOCsp}{\mathsf{LOC_{sp}}}
\newcommand{\sfLLDsp}{\mathsf{LLD_{sp}}}
\newcommand{\sfForget}{\mathsf{Forget}}
\newcommand{\sfnon}{\mathsf{non}}
\newcommand{\sfcov}{\mathsf{cov}}

\newcommand{\balpha}{\boldsymbol{\alpha}}
\newcommand{\biota}{\boldsymbol{\iota}}

\newcommand{\bB}{\mathbf{B}}
\newcommand{\bL}{\mathbf{L}}
\newcommand{\bg}{\mathbf{g}}
\newcommand{\bp}{\mathbf{p}}
\newcommand{\bq}{\mathbf{q}}
\newcommand{\br}{\mathbf{r}}

\DeclareMathOperator{\rmcard}{\mathrm{card}} 
\newcommand{\rmccc}{\mathrm{ccc}}
\DeclareMathOperator{\rmid}{\mathrm{id}} 
\DeclareMathOperator{\rmloc}{\mathrm{loc}}

\DeclareMathOperator{\rmTan}{\mathrm{Tan}} 

\newcommand{\hel} {
\hskip2.5pt{\vrule height7pt width.5pt depth0pt}
\hskip-.2pt\vbox{\hrule height.5pt width7pt depth0pt}
\, }

\newcommand{\oh}{\frac{1}{2}}

\renewcommand{\em}{\bf}
\renewcommand{\leq}{\leqslant}
\renewcommand{\geq}{\geqslant}
\renewcommand{\subset}{\subseteq}

\begin{document}


\title{A brief review of Radon-Nikod\'ymification}

\author[Th. De Pauw]{Thierry De Pauw}
\address{Institute for Theoretical Sciences,
Westlake University,
No. 600, Dunyu Road, Xihu District, Hangzhou, Zhejiang, 310030, China}

\email{thierry.depauw@westlake.edu.cn}

\keywords{Radon-Nikod\'ym theorem, Hausdorff measure, integral geometry, functions of bounded variation}

\subjclass[2010]{Primary 49Q15,28A75,26A45,28C05}




\maketitle


\tableofcontents


\section{Foreword}

In section \ref{sec.mot} we explain how a certain point of view on the dual of $\mathbf{SBV}(\Omega)$ is related to understanding the dual of $\mathbf{L}_1(\Omega, \calB(\Omega), \calH^{n-1}; \R^n)$, where $\calH^{n-1}$ is a Hausdorff measure.
Difficulties with the latter arise from the fact that the Radon-Nikod\'ym Theorem does not apply to the measure space $(\Omega, \calB(\Omega), \calH^{n-1})$.
In section \ref{sec.inj.sur} we review our paper \cite{DEP.21} discussing properties of the canonical map $\Upsilon : \bL_\infty(X,\calA,\mu) \to \bL_1(X,\calA,\mu)^*$ for a general measure space $(X,\calA,\mu)$: $\Upsilon$ is injective iff $(X,\calA,\mu)$ is semifinite; $\Upsilon$ is surjective iff $(X,\calA,\mu)$ is semilocalizable, {\it i.e.}\ the Boolean algebra $\frac{\calA}{\calN_{\mu,\rmloc}}$ is order-complete, where $\calN_{\mu,\rmloc}$ consists of those members of $\calA$ whose trace on each set of finite measure is null.
Injectivity is related to the uniqueness of Radon-Nikod\'ym derivatives whereas surjectivity is related to their existence.
We explain how, for Hausdorff measures, whether $\Upsilon$ is injective or surjective depends on the $\sigma$-algebra considered -- Borel or Hausdorff-measurable in the sense of Caratheodory.
In particular, we give hints as to why the surjectivity of $\Upsilon$ for the measure space $(\R^2,\calA_{\calH^1},\calH^1)$ is undecidable in ZFC.
In section \ref{sec.msn} we address the existence, proved in our paper \cite{BOU.DEP.21} with {\sc Ph. Bouafia}, for {\it every} measure space $(X,\calA,\mu)$ of a {\it Radon-Nikod\'ymification} $(\hat{X},\hat{\calA},\hat{\mu})$ akin the completion of a metric space or the compactification of a Hausdorff topological space. 
Rendering a measure space semifinite can (easily) be done without altering the base space $(X,\calA)$.
However, examples in section \ref{sec.inj.sur} show that rendering a measure space semilocalizable requires in general to replace $X$ by a ``larger'' space $\hat{X}$ -- a quotient of a certain fiber bundle.
In order to state the universal property defining the Radon-Nikod\'ymification of a measure space, one needs to set up the proper categorical language.
The general existence theorem is non-constructive.
In section \ref{sec.igm}, as in \cite{BOU.DEP.21} we give an explicit ``hands on'' description of $(\hat{X},\hat{\calA},\hat{\mu})$ in case $(X,\calA,\mu)$ is, for instance, an integral geometric measure space $(\R^2,\calB(\R^2),\calI^1_\infty)$.
This relies on specific properties of the measure $\calI^1_\infty$ (Besicovitch-Federer Projection Theorem, density 1) not shared by the Hausdorff measure $\calH^1$.
The replacement does not affect the application to the original problem of describing the dual of $\mathbf{SBV}(\Omega)$, since the total variation $\|Du\|$ of the distributional derivative $Du$ of a special function of bounded variation $u$ is absolutely continuous with respect to (the semifinite version of) $\calI^1_\infty$, which amounts to saying that the jump set of $u$ is countably rectifiable.
This leads, in section \ref{sec.sbv} to an integral representation formula for the members of $\mathbf{SBV}(\Omega)$ by means of vector fields defined on the enlarged space $\hat{\Omega}$, the main result of our paper \cite{DEP.BOU.22} with {\sc Ph. Bouafia}. 

\section{Motivation}
\label{sec.mot}

Let $n \geq 2$ and $\Omega \subset \Rn$ be open.
Some questions pertaining to the calculus of variations would benefit from a useful description of the dual of the Banach space $\mathbf{BV}(\Omega)$ of functions of bounded variation in the sense of {\sc E. De Giorgi} \cite[3.1]{AMBROSIO.FUSCO.PALLARA} and of the Banach space $\mathbf{SBV}(\Omega)$ of special functions of bounded variation \cite[4.1]{AMBROSIO.FUSCO.PALLARA}.
The first question occurs as Problem 7.4 in \cite{ARCATA} nearly forty years ago.
Measures belonging to the dual space of $\mathbf{BV}(\Omega)$ have been characterized by {\sc N.G. Meyers} and {\sc W.P. Ziemer} in \cite{MEY.ZIE.77} whereas signed measures belonging to this dual space have been characterized by {\sc N.C. Phuc} and {\sc M. Torres} \cite{PHU.TOR.17}.
A description of the other members was obtained (in a slightly different context) by {\sc F.J. Almgren} \cite{ALM.65.CH} under the Continuum Hypothesis and the particular description was proved to be independent of ZFC axioms by the present author in \cite{DEP.98}.
Recently, {\sc N. Fusco} and {\sc D. Spector} \cite{FUS.SPE.18} have given a more precise description under the Continuum Hypothesis, following former work of {\sc R.D. Mauldin}.
In the present note we review recent work \cite{DEP.21,BOU.DEP.21,DEP.BOU.22} that provides an explicit description of members of the dual of $\mathbf{SBV}(\Omega)$ within the ZFC axioms and that is ``optimal'' in some specific universal sense explained below.
\par 
For each $u$ representing a member of $\mathbf{SBV}(\Omega)$ the distributional gradient $Du$ of $u$ decomposes as $Du = \calL^n \hel \nabla u + \calH^{n-1} \hel j_u$ where
$\nabla u$ is the pointwise a.e. approximate gradient of $u$ and
$j_u$ is a vector field carried on the approximate discontinuity set
$S_u$ of $u$ on which we have $j_u = (u^+-u^-)\nu_u$, $\nu_u$ being a
unitary field normal to $S_u$ and $u^+, u^-$ the approximate limits
of $u$ on either sides of $S_u$.
We have
\begin{equation*}
  \|u\|_{\mathbf{SBV}(\Omega)} = \int_\Omega |u| \, d \calL^n + \int_{S_u} \|j_u\| \,
  d \calH^{n-1} + \int_\Omega \|\nabla u\| \, d \calL^n.
\end{equation*}
Thus, $u \mapsto (u, j_u, \nabla u)$ is an isometric embedding
\begin{equation*}
  \mathbf{SBV}(\Omega) \to \mathbf{L}_1(\Omega, \calB(\Omega), \calL^n) \times
  \mathbf{L}_1(\Omega, \calB(\Omega), \calH^{n-1}; \R^n) \times \mathbf{L}_1(\Omega,
  \calB(\Omega), \calL^n ; \R^n).
\end{equation*}
Since $(\Omega,\calB(\Omega),\calL^n)$ is $\sigma$-finite, the dual of the corresponding $\mathbf{L}_1$ space is, of course, the corresponding $\mathbf{L}_\infty$ space. 
This remark does not apply to the middle measure space $(\Omega,\calB(\Omega),\calH^{n-1})$. 
An attempt at understanding the dual of the corresponding $\mathbf{L}_1$ space leads to a study of the canonical map
\begin{equation*}
\Upsilon \colon \mathbf{L}_\infty(\Omega,\calB(\Omega),\calH^{n-1}) \to \mathbf{L}_1(\Omega,\calB(\Omega),\calH^{n-1})^*\,.
\end{equation*}

\section{Injectivity and surjectivity of the canonical map $\bL_\infty \to \bL_1^*$}
\label{sec.inj.sur}

Let $(X,\calA,\mu)$ be a measure space.
There is a natural linear retraction
\begin{equation}
\label{eq.100}
\Upsilon : \bL_\infty(X,\calA,\mu) \to \bL_1(X,\calA,\mu)^*
\end{equation}
which sends $\bg$ to $\mathbf{f} \mapsto \int_X gf d\mu$, where $g$ and $f$ represent $\bg$ and $\mathbf{f}$, respectively.
In general, $\Upsilon$ may not be injective nor surjective.
It has been understood for a long time that $\Upsilon$ is injective if and only if $(X,\calA,\mu)$ is {\it semifinite}.
This means that each $A \in \calA$ such that $\mu(A)=\infty$ admits a subset $\calA \ni B \subset A$ with $0 < \mu(B) < \infty$.
Of course, every $\sigma$-finite measure space is semifinite.
Yet, the dependence upon the $\sigma$-algebra under consideration already occurs in the case of interest to us.
The situation is the following.
\begin{enumerate}
\item If $X$ is a complete separable metric space and $0 < d < \infty$, then the measure space $(X,\calB(X),\calH^d)$ is semifinite. Here, $\calB(X)$ denotes the $\sigma$-algebra of Borel subsets of $X$ and $\calH^d$ is the $d$-dimensional Hausdorff measure on $X$. In case $X=\Rn$, this was proved by {\sc R.O. Davies} \cite{DAV.52} and in general by {\sc J. Howroyd} \cite{HOW.95}.
\item According to {\sc D.H. Fremlin} \cite[439H]{FREMLIN.IV}, the measure space $(\R^2,\calA_{\calH^1},\calH^1)$ is not semifinite, where $\calA_{\calH^1}$ denotes the $\sigma$-algebra consisting of $\calH^1$-measurable subsets of $\R^2$. This is based on the existence of ``large'' universally null subsets of $[0,1]$, established by {\sc E. Grzegorek}, \cite{GRZ.81}. See also the article of {\sc O. Zindulka} \cite{ZIN.12}.
\end{enumerate}
\par
Under the assumption that $(X,\calA,\mu)$ is semifinite, a necessary and sufficient condition for the surjectivity of $\Upsilon$ has been known for a long time.
It asks for the quotient Boolean algebra $\calA/\calN_\mu$ to be order complete, where $\calN_\mu = \calA \cap \{ N : \mu(N) = 0 \}$ is the $\sigma$-ideal of $\mu$-null sets.
Semifinite measure spaces with this property are sometimes called {\it Maharam}, \cite[211G]{FREMLIN.II}.
A stronger condition, sometimes called {\it decomposable}, generalizes the idea of $\sigma$-finiteness to possibly uncountable decomposition into sets of finite measure, together with a new condition called {\it locally determined} (that measurability be determined by sets of finite measure), see \cite[211E]{FREMLIN.II} for the definition of decomposable.
If the quotient $\sigma$-algebra $\calA/\calN_\mu$ is not too big, then decomposability implies Maharam, according to {\sc E.J. McShane} \cite{MCS.62}, but not in general, according to {\sc D.H. Fremlin} \cite[216E]{FREMLIN.II}.
\par 
If $X$ is a Polish space and $\calB(X)$ denotes the $\sigma$-algebra consisting of its Borel subsets, and if the measure space $(X,\calB(X),\mu)$ is decomposable, then it is $\sigma$-finite. 
I learned the ``counting argument'' to prove this from {\sc D.H. Fremlin}, see \cite{DEP.21}.
In view of (1) above, this shows that $(\R^2,\calB(\R^2),\calH^1)$ is not decomposable.
Since decomposability is stronger in general than the surjectivity of $\Upsilon$, we need to argue a bit more to show that $(\R^2,\calB(\R^2),\calH^1)$ is not Maharam, see below.
This observation calls for developing a criterion for the surjectivity of $\Upsilon$, without assuming that $(X,\calA,\mu)$ be semifinite in the first place.
%
%
Regarding the question whether 
\begin{equation*}
\Upsilon : \bL_\infty\left(\R^2,\calA,\calH^1\right) \to \bL_1\left(\R^2,\calA,\calH^1\right)^*
\end{equation*}
is surjective or not, the situation is the following, see \cite{DEP.21}.
\begin{enumerate}
\item[(3)] If $\calA=\calB(\R^2)$, then $\Upsilon$ is not surjective. Since $(\R^2,\calB(\R^2),\calH^1)$ is semifinite, according to (1), and not $\sigma$-finite, it is not decomposable. The argument of {\sc E.J. McShane} does not show that $(\R^2,\calB(\R^2),\calH^1)$ is not Maharam (the reason being that its completion is not locally determined). However, we give below a simple argument to the extent that it is not Maharam, based on Fubini's Theorem.
\item[(4)] If $\calA=\calA_{\calH^1}$, then the surjectivity of $\Upsilon$ is undecidable in $\mathsf{ZFC}$. The consistency of its surjectivity is a consequence of the Continuum Hypothesis. The consistency of it not being surjective was first noted in \cite{DEP.98}, although in a slightly different disguise. The idea is explained below.
\end{enumerate}
Why, however, would the answer depend upon the $\sigma$-algebra under consideration?
In order to understand this, let us try to prove that $\Upsilon$ is surjective.
\par 
We know from the classical Riesz' Theorem that $\Upsilon$ is surjective whenever $(X,\calA,\mu)$ is a finite measure space -- and it boils down to applying the Radon-Nikod\'ym Theorem.
This suggests to consider $\calA^f_\mu = \calA \cap \{ A : \mu(A) < \infty \}$ and, for each $A \in \calA^f_\mu$, the map
\begin{equation*}
\Upsilon^A : \bL_\infty(A,\calA_A,\mu_A) \to \bL_1(A,\calA_A,\mu_A)^*,
\end{equation*}
where $(A,\calA_A,\mu_A)$ is the obvious (finite) measure subspace.
Thus, $\Upsilon^A$ is an isometric linear isomorphism and, given $\alpha \in  \bL_1(X,\calA,\mu)^*$, there exist $g_A \in \bg_A \in \bL_\infty(A,\calA_A,\mu_A)$ such that
\begin{equation*}
(\alpha \circ \iota_A)(\mathbf{f}) = \int_A g_Af d\mu_A,
\end{equation*}
whenever $f \in \mathbf{f} \in \bL_1(A,\calA_A,\mu_A)$, where $\iota_A : \bL_1(A,\calA_A,\mu_A) \to \bL_1(X,\calA,\mu)$ is the obvious embedding.
From the $\mu_A$-almost everywhere uniqueness of the Radon-Nikod\'ym derivative $g_A$, we infer that if $A,A' \in \calA^f_\mu$, then $\mu ( A \cap A' \cap \{ g_A \neq g_{A'} \}) = 0 $.
Thus, $(g_A)_{A \in \calA^f_\mu}$ is what we call, from now on, a {\it compatible family} of locally defined measurable functions and the question is whether it corresponds to a globally defined measurable function, {\it i.e.}\ whether there exists an $\calA$-measurable $g : X \to \R$ such that $\mu(A \cap \{ g \neq g_A \}) = 0$, for every $A \in \calA^f_\mu$.
If such $g$ exists, let us call it a {\it gluing} of the compatible family $(g_A)_{A \in \calA^f_\mu}$.
The problem of gluing a compatible family in our setting is
reminiscent of the fact that, for a topological space $X$, the functor
of continuous functions on open sets is a sheaf.
However, unlike in the case of continuous functions, in order to
define $g$ globally, we ought to make choices on the domains $A \cap
A'$, for $A,A' \in \calA^f_\mu$, because $g_A$ and $g_{A'}$ do not coincide
everywhere there, but merely {\it almost everywhere}.
In general, one cannot make consistent choices. 
\par 
It turns out to be rather useful to notice that the question whether a gluing exists or not can be asked in a slightly more general setting since it depends on the measure $\mu$ only insofar as its $\mu$-null sets are involved.
Thus, a {\it measurable space with negligibles} $(X,\calA,\calN)$, abbreviated as MSN, consists of a measurable space $(X,\calA)$ and a $\sigma$-ideal $\calN \subset \calA$.
Given any $\calE \subset \calA$, one can readily define the notion of a compatible family $(g_E)_{E \in \calE}$ by asking that $E \cap E' \cap \{ g_E \neq g_{E'} \} \in \calN$ whenever $E,E' \in \calE$, and by saying that an $\calA$-measurable function $g : X \to \R$ is a gluing of $(g_E)_{E \in \calE}$ provided $E \cap \{ g \neq g_E \} \in  \calN$, for all $E \in \calE$.
One then shows that each compatible family admits a gluing if and only if each $\calE \subset \calA$ admits an $\calN$-essential supremum $A \in \calA$.
This means that 
\begin{enumerate}
\item[(i)] For every $E \in \calE$ one has $E \setminus A \in \calN$;
\item[(ii)] For every $B \in \calA$, if $E \setminus B \in \calN$ whenever $E \in \calE$, then $A \setminus B \in \calN$. 
\end{enumerate}
We say that a measurable space with negligibles is {\it localizable} if it has this property.
\par 
In \cite{DEP.21} we characterize those measure spaces such that $\Upsilon$ is surjective.
To state this, we first define
\begin{equation*}
\calN_\mu\left[\calA^f_\mu\right] = \calA \cap \left\{ N : \mu(A \cap N) = 0 \text{ for all } A \in \calA^f_\mu \right\} .
\end{equation*}
It is a $\sigma$-ideal, whose members one is tempted to call {\it locally $\mu$-null}.
We will also denote it simply as $\calN_{\mu,\rmloc}$.
\begin{Theorem*}
For any measure space $(X,\calA,\mu)$, the map $\Upsilon$ (recall \eqref{eq.100}) is surjective if and only if the measurable space with negligibles $\left(X,\calA, \calN_{\mu,\rmloc}\right)$ is localizable.
\end{Theorem*}
We call a measure space {\it semilocalizable} if it has this property -- no semifiniteness is assumed.
In \cite{DEP.21} we study the connection with the notion of {\it almost decomposable} measure space introduced in \cite{DEP.98} thereby generalizing to non semifinite measure spaces the classical theory briefly evoked above.
We call a measure space $(X,\calA,\mu)$ {\it almost decomposable} if there exists a disjointed family $\calG \subset \calA^f_\mu$ such that
\begin{equation*}
\forall A \in \calP(X) : ( \forall G \in \calG : A \cap G \in \calA) \Rightarrow A \in \calA \,,
\end{equation*}
and
\begin{equation*}
\forall A \in \calA : \mu(A) < \infty \Rightarrow \mu(A) = \sum_{G \in \calG} \mu(A \cap G) \,.
\end{equation*}
Using an idea of {\sc E.J. McShane} \cite{MCS.62} and the fact that there are not too many equivalence classes of measurable sets with respect to a Borel regular outer measure on a Polish space, we prove the following \cite{DEP.21}.
\begin{Theorem*}
Let $X$ be a complete separable metric space and $0 < d < \infty$. For the measure space $(X,\calA_{\calH^d},\calH^d)$ the following are equivalent.
\begin{enumerate}
\item The canonical map $\Upsilon$ is surjective;
\item $(X,\calA_{\calH^d},\calH^d)$ is semilocalizable;
\item $(X,\calA_{\calH^d},\calH^d)$ is almost decomposable.
\end{enumerate}
\end{Theorem*}
\par 
Let us now consider the measure space $(\R^2,\calB(\R^2),\calH^1)$ in view of the notion of semilocalizability.
We know it is not semilocalizable, (3) above, and we promised to show how this is a consequence of Fubini's Theorem.
Define the vertical sections $V_s = \{s\} \times \R$, $s \in \R$, and the horizontal sections $H_t = \R \times \{t\}$, $t \in \R$. 
Assume if possible that $A \in \calB(\R^2)$ is an $\calN_{\calH^1}\left[\calB(\R^2)^f_{\calH^1}\right]$-essential supremum of the family $(V_s)_{s \in \R}$.
It would then readily follow that
\begin{enumerate}
\item[(a)] $\calH^1(V_s \setminus A) = 0$ for every $s \in \R$;
\item[(b)] $\calH^1(H_t \cap A) = 0$ for every $t \in \R$.
\end{enumerate}
Indeed, upon noticing that $V_s$ and $H_t$ have $\sigma$-finite $\calH^1$ measure, (a) is a rephrasing of (i) above and (b) follows from (ii) applied with $B=A \setminus H_t$.
Applying Fubini's Theorem twice would yield
\begin{equation*}
\calL^2(\R^2 \setminus A) = \int_\R \calH^1( V_s \setminus A) d\calL^1(s) = 0,
\end{equation*}
according to (a), and
\begin{equation*}
\calL^2(\R^2 \cap A) = \int_\R \calH^1( H_t \cap A) d\calL^1(t) = 0,
\end{equation*}
according to (b). 
In turn, $\calL^2(\R^2)=0$, a contradiction. 
Clearly, the same argument applies with $\R^2$ replaced by any Borel set $X \subset \R^2$ such that $\calL^2(X) > 0$, to showing that $(X,\calB(X),\calH^1)$ is not semilocalizable.
\par 
There are two cases when the above argument is not conclusive: 
\begin{enumerate}
\item[($\alpha$)] when $A$ is not $\calL^2$-measurable (because Fubini's Theorem does not apply);
\item[($\beta$)] when $\calL^2(X)=0$ (because no contradiction ensues).
\end{enumerate}
\par 
With regard to case ($\alpha$), indeed, when we replace the $\sigma$-algebra $\calB(\R^2)$ by the larger $\calA_{\calH^1}$, then $(\R^2,\calA_{\calH^1},\calH^1)$ is consistently semilocalizable.
This is a consequence of the Continuum Hypothesis, details are in \cite{DEP.21}; in fact, a much more general statement holds\footnote{I learned it from \cite[2.5.10]{GMT}. Unfortunately the presentation there does not allow for putting emphasis on the role played by the choice of a particular $\sigma$-algebra.}.
As noticed in \cite{DEP.98}, it turns out, however, that $(\R^2 , \calA_{\calH^1}, \calH^1)$ is {\it also} consistently {\it not} semilocalizable.
Here is the reason why.
We assume that $A \in \calA_{\calH^1}$ is an $\calN_{\calH^1}\left[\calA_{\calH^1}^f\right]$-essential supremum of the family $(V_s)_{s \in \R}$.
For each $s \in \R$, we define $T_s = \R \cap \{ t : (s,t) \in V_s \setminus A \}$, thus, $\calL^1(T_s)=0$, according to (a).
Now, choose $E \subset \R$ such that $\calL^1(E) > 0$ and $E$ has least cardinal among all sets with nonzero Lebesgue measure, and let $\sfnon(\calN_{\calL^1})$ denote this cardinal.
Assume that there exists $t \in \R \setminus \cup_{s \in E} T_s$.
Then, for each $s \in E$, $t \not \in T_s$, {\it i.e.}\ $(s,t) \in H_t \cap A$.
Therefore, $\calL^1(E) = 0$, according to (b), a contradiction.
Of course, we can reach this contradiction only if $\R \neq \cup_{s \in E} T_s$, which depends upon how big $E$ is.
We denote by $\sfcov(\calN_{\calL^1})$ the least cardinal of a covering of $\R$ by $\calL^1$-negligible sets. 
Thus, if $\rmcard E = \sfnon(\calN_{\calL^1}) < \sfcov(\calN_{\calL^1})$, then the argument goes through.
It turns out that this strict inequality of cardinals (appearing in the so-called Cicho\'n diagram) is consistent with $\mathsf{ZFC}$, \cite[Chapter 7]{BARTOSZYNSKI.JUDAH} or \cite[552H and 552G]{FREMLIN.V.2}.
We will refer to this idea below as the ``vertical-horizontal method''.
This argument is from \cite{DEP.98}; I learned it from {\sc D.H. Fremlin}.
\par 
With regard to case ($\beta$) above, we observe again that the Continuum Hypothesis implies that $(X,\calA_{\calH^1},\calH^1)$ is semilocalizable, for any compact set $X \subset \R^2$, regardless whether it has zero $\calL^2$ measure or not.
The question is, therefore, whether $(X,\calA_{\calH^1},\calH^1)$ is semilocalizable in $\mathsf{ZFC}$ or consistently not semilocalizable.
The latter occurs when the ``vertical-horizontal method'' generalizes from $X = \R^2$ to $X$.
For instance, it clearly generalizes to $X = [a,b ] \times [c,d]$, but it is not instantly obvious how to proceed if $\calL^2(X) = 0$.
%
Thus, we ought to explain how the ``vertical-horizontal method'' described above, showing that if $\sfnon(\calN_{\calL^1}) < \sfcov(\calN_{\calL^1})$ then $\left(\R^2,\calA_{\calH^1},\calH^1\right)$ is not semilocalizable, can be adapted to the case where $\R^2$ is replaced with some suitable subset $X \subset \R^2$.
First of all, it is useful to observe that if $(S,\calB(S),\sigma)$ is a probability space, $S$ is Polish, and $\sigma$ is diffuse, then $\sfnon(\calN_{\bar{\sigma}}) = \sfnon(\calN_{\calL^1})$ and $\sfcov(\calN_{\bar{\sigma}}) = \sfcov(\calN_{\calL^1})$.
This ensues from the Kuratowski Isomorphism Theorem \cite[theorem 3.4.23]{SRIVASTAVA}.
A careful inspection of the argument leads to the following see \cite{DEP.21}.
\begin{Theorem*}
Let $0 < d < 1$ and let $C_d \subset [0,1]$ be the standard self-similar Cantor set of Hausdorff dimension $0 < d < 1$.
Whether the measure space $(C_d \times C_d , \calA_{\calH^d}, \calH^d)$ is semilocalizable is undecidable in $\mathsf{ZFC}$.
\end{Theorem*}
Incidentally, constructing a certain isomorphism in the category of measurable spaces with negligibles, we are able to infer the following as well.
\begin{Theorem*}
Whether the measure space $\left([0,1],\calA_{\calH^{1/2}},\calH^\oh\right)$ is semilocalizable is undecidable in $\mathsf{ZFC}$.
\end{Theorem*}
Here, the exponent $1/2$ reflects the nature of the argument, viewing the space $X$ as a product of a kind, where ``vertical'' sets $V_s$ and ``horizontal'' sets $H_t$ of the same size make sense.
Note that the sets $X = C_d \times C_d$ are purely $(\calH^1,1)$-unrectifiable.
We next seek to apply the ``vertical-horizontal method'' to an $\calL^2$-negligible compact set $X \subset \R^2$ which is not purely $(\calH^1,1)$-unrectifiable and prove that $(X,\calA_{\calH^1},\calH^1)$ is not semilocalizable.
Let us choose $X$ as small as possible, {\it i.e.}\ of Hausdorff dimension 1, say $X = C \times [0,1]$, where $C \subset [0,1]$ is a Cantor set of Hausdorff dimension 0.
It is, of course, clear that $V_s = \{s\} \times [0,1]$, $s \in C$, can be chosen as our vertical sets, yet the choice $H_t = C \times \{t\}$, $t \in [0,1]$, will be of no use since $\calH^1(H_t)=0$ and, therefore, no contradiction can ensue when implementing the ``vertical-horizontal method''.
Instead, we proceed as follows to define $H_t$.
Let $\mu$ be a diffuse probability measure on $C$ and let $f(t) = \mu([0,t])$, $t \in [0,1]$, be its distribution function (this is a version of the Devil's staircase corresponding to our 0 dimensional set $C$).
Consider the graph $G$ of the function $\frac{1}{2}f$ ; thus $G$ is a rectifiable curve, and intersects non $\calH^1$ trivially the set $X$.
We then define $H_t = G + t.e_2$, $t \in [0,1/2]$, where $e_2=(0,1)$. 
It turns out that these will successfully play the role of horizontal sets.
We obtain the following \cite{DEP.21}.
\begin{Theorem*}
Assume that
\begin{enumerate}
\item $C \subset [0,1]$ is some Cantor set of Hausdorff dimension 0; 
\item $X = C \times [0,1]$;
\item $\calA$ is a $\sigma$-algebra and $\calB(X) \subset \calA \subset \calP(X)$;
\item $\calN = \calN_{\calH^1}$ or $\calN = \calN_{pu}$. 
\end{enumerate}
It follows that the measurable space with negligibles $(X,\calA,\calN)$ is consistently not localizable.
\end{Theorem*}
Here, $\calN_{pu}$ is the $\sigma$-ideal consisting of those subsets $S$ of $\R^2$ that are purely $(\calH^1,1)$-unrectifiable, {\it i.e.}\ $\calH^1(S \cap \Gamma) = 0$ for every
Lipschitz (or, for that matter, $C^1$) curve $\Gamma \subset \R^2$.

\section{Measurable spaces with negligibles}
\label{sec.msn}

It is natural to want to associate, with an arbitrary measure space
$(X,\calA,\mu)$, an improved version of itself -- in a universal way
-- ideally one for which the Radon-Nikod\'ym Theorem holds.
According to the previous section, this amounts to making it become semifinite and semilocalizable.
It is not difficult to modify slightly the measure $\mu$, keeping the
underlying measurable space $(X,\calA)$ untouched, in order to make it
semifinite.
Specifically, letting $\mu_{\mathrm{sf}}(A) = \sup \{ \mu(A \cap F) :
F \in \calA^f \}$, for $A \in \calA$, one checks that
$(X,\calA,\mu_{\mathrm{sf}})$ is semifinite and that
$\calN_{\mu_{\mathrm{sf}}} = \calN_{\mu,\rmloc}$.
However, it appears to be a more delicate task to modify
$(X,\calA,\mu)$ in a canonical way in order for $\Upsilon$ to become
surjective.
One may naively
attempt to ``add measurable sets'' in a smart way in order to obtain a
localizable MSN $(X,\hat{\calA},\hat{\calN}_{\hat{\mu}})$, just as many as needed,
and that would be a ``localizable version'' of $(X,\calA,\calN_\mu)$.
Unfortunately, within ZFC this cannot always be done while ``sticking
in the base space $X$'', as shown by the last theorem stated in the previous section.
\par 
Thus, one may need to also add points to the base space $X$.
In the case of a general measure space $(X,\calA,\mu)$, we can get a
feeling of what needs to be done, when trying to define the gluing of
a compatible family $( f_F )_{F \in \calA^f}$.
Indeed, each $x \in X$ may belong to several $F \in \calA^f$ and this
calls for considering an appropriate quotient of the fiber bundle
$\{(x,F) : x \in F \in \calA^f\}$.  
\par 
Our task is to define a general notion of
``localization'' of an MSN and to prove existence results in some
cases.
This will be akin the ``completion'' of a metric space and the ``compactification'' of a Hausdorff topological space.
Both cases can be described ``by hand'' in various ways but exist only up to an isomorphism in a relevant category and the fact that one adds ``as few elements as needed'' in each case is conveniently expressed by a universal property.
Since a definition of ``localization'' of an MSN will also involve a universal
property, it is critical to determine which category is appropriate
for our purposes.
As this offers unexpected surprises, we describe the several steps in
some detail.
\par 
The objects of our first category $\sfMSN$ are the {\em saturated} MSN
$(X,\calA,\calN)$, by what we mean that for every $N,N' \subset X$, if
$N \subset N'$ and $N' \in \calN$, then $N \in \calN$.
This is in analogy with the notion of a complete measure space. 
In order to define the morphisms between two objects $(X,\calA,\calN)$
and $(Y,\calB,\calM)$, we say that a map $f \colon X \to Y$ is
$[(\calA,\calN),(\calB,\calM)]$-measurable if $f^{-1}(B) \in \calA$
for every $B \in \calB$ and $f^{-1}(M) \in \calN$ for every $M \in
\calM$.
For instance, if $X$ is a Polish space and $\mu$ is a diffuse
probability measure on $X$, there exists (by the Kuratowski Isomorphism Theorem) a
Borel isomorphism $f \colon X \to [0,1]$ such that $f_\# \mu = \calL^1$,
where $\calL^1$ is the Lebesgue measure, thus $f$ is
$[(\calB(X),\calN_\mu),(\calB([0,1]),\calN_{\calL^1})]$-measurable.
We define an equivalence relation for such measurable maps $f,f' \colon X
\to Y$ by saying that $f \sim f'$ if and only if $\{ f \neq f' \} \in
\calN$.
The morphisms in the category $\sfMSN$ between the objects
$(X,\calA,\calN)$ and $(Y,\calB,\calM)$ are the equivalence classes of
$[(\calA,\calN),(\calB,\calM)]$-measurable maps.
At this stage, we need to suppose that $(X, \calA, \calN)$ is
saturated for the relation of equality almost everywhere to be
transitive. With this assumption, the composition of
measurable maps is also compatible with $\sim$.
\par We let $\sfLOC$ be the full subcategory of $\sfMSN$ whose objects
are the localizable MSNs.
We may be tempted to define the localization of a saturated MSN
$(X,\calA,\calN)$ as its coreflection (if it exists) along the
forgetful functor $\sfForget \colon \sfLOC \to \sfMSN$, and the question of
existence in general becomes that of the existence of a right adjoint
to $\sfForget$.
Specifically, we may want to say that a pair
$[(\hat{X},\hat{\calA},\hat{\calN}),\bp]$, where
$(\hat{X},\hat{\calA},\hat{\calN})$ is saturated localizable MSN and
$\bp$ is a morphism $\hat{X} \to X$, is a localization of
$(X,\calA,\calN)$ whenever the following universal property holds.
For every pair $[(Y,\calB,\calM),\bq]$, where $(Y,\calB,\calM)$ is a
saturated localizable MSN and $\bq$ is a morphism $Y \to X$, there
exists a unique morphism $\br \colon Y \to \hat{X}$ such that $\bq = \bp
\circ \br$.
\begin{equation}
\label{diagram}
  \begin{tikzcd}
    (Y, \calB, \calM) \arrow[rd, swap, "\bq"] \arrow[rr, dotted,
      "\exists!\br"] & & (\hat{X}, \hat{\calA}, \hat{\calN})
    \arrow[ld, "\bp"] \\ & (X, \calA, \calN) &
  \end{tikzcd} \tag{$\bigtriangledown$}
\end{equation}
\par However, we now illustrate that the notion of morphism defined so
far is not yet the appropriate one that we are after.
We consider the MSN $(X,\calA,\{\emptyset\})$ where $X = \R$ and
$\calA$ is the $\sigma$-algebra of Lebesgue measurable subsets of
$\R$.
We observe that the only reasonable choice for the localization of $(X,\calA,\{\emptyset\})$
is $[(X,\calP(X),\{\emptyset\}),\bp]$ with $\bp$ induced by the
identity $\rmid_X$.
Assume if possible that this is the case.
In the diagram above we consider $(Y, \calB,\calM) =
(X,\calA,\calN_{\calL^1})$ and $\bq$ induced by the identity.
Note that this is, indeed, a localizable MSN since it is associated
with a $\sigma$-finite measure space (this is folklore and can be found for instance in \cite{DEP.21} and \cite{DEP.BOU.22}).
Thus, there would exist a morphism $\br$ in $\sfMSN$ such that $\bp
\circ \br = \bq$.
Picking $r \in \br$, this implies that $X \cap \{ x: r(x) \neq x \}$
is Lebesgue negligible.
The measurability of $r$ would then imply that $r^{-1}(S) \in \calA$
for every $S \in \calP(X)$, contradicting the existence of non
Lebesgue measurable subsets of $\R$.
\par The problem with the example above is that the objects
$(X,\calA,\calN_{\calL^1})$ and $(X,\calA,\{\emptyset\})$ should not
be compared, in other words that $\bq$ should not be a legitimate
morphism.
We say that a morphism $\mathbf{f} \colon (X,\calA,\calN) \to
(Y,\calB,\calM)$ of the category $\sfMSN$ is {\em supremum preserving}
if the following holds for (one and, therefore, every) $f
\in\mathbf{f}$.
If $\calF \subset \calB$ admits an $\calM$-essential supremum $S \in
\calB$, then $f^{-1}(S)$ is an $\calN$-essential supremum of
$f^{-1}(\calF)$.
It is easy to see that adding this condition to the definition of
morphism rules out the $\bq$ considered in the preceding paragraph.
We define the category $\sfMSNsp$ to be that whose objects are the
saturated MSNs and whose morphisms are those morphisms of $\sfMSN$
that are supremum preserving.
We define similarly $\sfLOCsp$.
We now define the {\em localizable version} (if it exists) of a
saturated MSN with the similar universal property illustrated in
\eqref{diagram}, except for we now require all morphisms to be in
$\sfMSNsp$, {\it i.e.}\ supremum preserving.
In other words, it is a coreflection of an object of $\sfMSNsp$ along
$\sfForget \colon \sfLOCsp \to \sfMSNsp$.
Unfortunately, this is not quite yet the right setting.
Indeed, we show in \cite{DEP.BOU.22} that if $X$ is uncountable
and $\calC(X)$ is the countable-cocountable $\sigma$-algebra of $X$,
then $[(X,\calP(X),\{\emptyset\}),\biota]$ (with $\biota$ induced by
$\rmid_X$) is not the localizable version of
$(X,\calC(X),\{\emptyset\})$.
This prompts us to introduce a new category.
\par 
In an MSN $(X,\calA,\calN)$ we say that a collection $\calE \subset \calA$ is {\em $\calN$-generating } if it admits $X$ as an $\calN$-essential supremum.
We say that an object $(X,\calA,\calN)$ of $\sfMSN$ is {\em
  locally determined} if for every $\calN$-generating collection
$\calE \subset \calA$ the following holds:
\begin{equation*}
  \forall A \subset X : \big[ \forall E \in \calE : A \cap E \in \calA
    \big] \Rightarrow A \in \calA.
\end{equation*}
We say that an object of $\sfMSN$ is {\em lld} if it is both
localizable and locally determined, and we let $\sfLLDsp$ be the
corresponding full subcategory of $\sfLOCsp$.
We now define the {\em lld version} of an object of $\sfMSNsp$ to be
its coreflection (if it exists) along $\sfForget \colon \sfLLDsp \to
\sfMSNsp$, {\it i.e.}\ it satisfies the corresponding universal property
illustrated in \eqref{diagram} with $Y$ and $\hat{X}$ being lld, and
the morphisms being supremum preserving.
This definition is satisfactory in at least the simplest case: If $(X,\calA,\{\emptyset\})$ is so that $\calA$
contains all singletons, then it admits
$[(X,\calP(X),\{\emptyset\}),\biota]$ as its lld version.
\par Our general question has now become whether $\sfForget : \sfLLDsp
\to \sfMSNsp$ admits a right adjoint.
Freyd's Adjoint Functor Theorem \cite[3.3.3]{BORCEUX.1} could prove
useful, however, we do not know whether it applies, mostly because we do
not know whether coequalizers (quotients) exist in $\sfMSNsp$.
We gather in Table \ref{table} the information that we know about
limits and colimits in the three categories we introduced.

\begin{table}[h]
\label{table}
\begin{tabular}{|c|c|c|c|}
\hline
& $\sfMSNsp$ & $\sfLOCsp$ & $\sfLLDsp$  \\ \hline
\hline
equalizers &  exist if $\{f=g\}$ is meas. (C) & $\boldsymbol{?}$ & exist \\ \hline
products & (countable) exist  & $\boldsymbol{?}$ & $\boldsymbol{?}$ \\ \hline
coequalizers & $\boldsymbol{?}$ & $\boldsymbol{?}$  & $\boldsymbol{?}$ \\ \hline
coproducts & exist (D) & exist  & exist (D)\\ \hline
\end{tabular}\hskip 1cm
\vskip .3cm
\caption{Limits and colimits in the three categories of MSNs.}
\end{table}
\par In view of proving some partial existence result for lld
versions, we introduce the intermediary notion of a {\em cccc}
saturated MSN, short for coproduct (in $\sfMSNsp$) of ccc saturated
MSNs.
It is easy to see that cccc MSNs are lld.
The {\em cccc version} of an object of $\sfMSNsp$ is likewise defined
by its universal property in diagram \eqref{diagram}, using supremum
preserving morphisms.
Our main results are about {\em locally ccc} MSNs, {\it i.e.}\ those
saturated MSNs $(X,\calA,\calN)$ such that $\calE_{\rmccc} = \calA
\cap \{ Z : \text{ the subMSN } (Z,\calA_Z,\calN_Z) \text{ is ccc}\}$
is $\calN$-generating.
A complete semifinite measure space $(X,\calA,\mu)$ is clearly
locally ccc, since $\calA^f$ is $\calN_\mu$-generating.
\begin{Theorem*}
  Let $(X,\calA,\calN)$ be a saturated locally ccc MSN. The following
  hold.
  \begin{enumerate}
  \item[(1)] $(X,\calA,\calN)$ admits a cccc version.
  \item[(2)] If furthermore $\calE_{\rmccc}$ contains an
    $\calN$-generating subcollection $\calE$ such that $\rmcard \calE
    \leq \mathfrak{c}$ and each $(Z,\calA_Z)$ is countably separated,
    for $Z \in \calE$, then $(X,\calA,\calN)$ admits an lld version
    which is also its cccc version.
  \end{enumerate}
\end{Theorem*}
\par By saying that a measurable space $(Z,\calA_Z)$ is countably
separated we mean that $\calA_Z$ contains a countable subcollection
that separates points in $Z$.
The cccc version $(\hat{X},\hat{\calA},\hat{\calN})$ is obtained as a
coproduct $\coprod_{Z \in \calE} (Z,\calA_Z,\calN_Z)$ where $\calE$ is
an $\calN$-generating almost disjointed refinement of
$\calE_{\rmccc}$, whose existence ensues from Zorn's Lemma.
In order to establish that this, in fact, is also the lld version
under the extra assumptions in (2), we need to build an appropriate
morphism $\br$ in diagram \eqref{diagram}, associated with an lld pair
$[(Y,\calB,\calM),\bq]$.
It is obtained as a gluing of $( q_Z )_{Z \in \calE}$ where $q_Z
\colon q^{-1}(Z) \to \hat{X}$ is the obvious map.
Since $\calE$ is almost disjointed, $( q_Z )_{Z \in \calE}$ is
compatible and, since $(Y,\calB,\calM)$ in diagram \eqref{diagram} is
localizable, the only obstruction to gluing is that $\hat{X}$ is not
$\R$.
Notwithstanding, $(\hat{X},\hat{\calA}) = \coprod_{Z \in \calE}
(Z,\calA_Z)$ is itself countably separated because $\rmcard \calE \leq
\mathfrak{c}$ so that the local determinacy of
$(Y,\calB,\calM)$ and the fact that $q^{-1}(\calE)$ is
$\calM$-generating (because $\calE$ is $\calN$-generating and $q$ is
supremum preserving) provides a gluing $r$.
\par We now explain how this applies to associating, in a canonical
way, a decomposable measure space \cite[211E]{FREMLIN.II} with {\it any} measure space
$(X,\calA,\mu)$ (the Radon-Nikod\'ym Theorem holds for decomposable measure spaces).
First, we recall that without changing the base space $X$ we can
render the measure space complete and semifinite.
In that case, $\calA^f$ is $\calN_\mu$-generating and witnesses the
fact that the saturated MSN $(X,\calA,\calN_\mu)$ is locally ccc.
By the theorem above, it admits a cccc version
$[(\hat{X},\hat{\calA},\hat{\calN}),\bp]$.
\begin{Theorem*}
  Let $(X,\calA,\mu)$ be a complete semifinite measure space and
  $[(\hat{X},\hat{\calA},\hat{\calN}),\bp]$ its corresponding cccc
  version. Let $p \in \bp$. There exists a unique (and independent of
  the choice of $p$) measure $\hat{\mu}$ defined on $\hat{\calA}$ such
  that $p_\# \hat{\mu} = \mu$ and $\calN_{\hat{\mu}}=\hat{\calN}$.
  Furthermore $(\hat{X},\hat{\calA},\hat{\mu})$ is a decomposable measure space, and the Banach spaces
  $\bL_1(X,\calA,\calN)$ and $\bL_1(\hat{X},\hat{\calA},\hat{\mu})$
  are isometrically isomorphic.
\end{Theorem*}
We call $(\hat{X},\hat{\calA},\hat{\mu})$ the {\em Radon-Nikod\'ymification} of $(X,\calA,\mu)$.

\section{Integral geometric measure}
\label{sec.igm}

Of course, the general process for constructing $\hat{X}$ in the last section is non-constructive, as it involves the axiom of choice to turn $\calA^f$
into an almost disjointed generating family.
This is why, here, we briefly explore a
particular case where we are able to describe explicitly $\hat{X}$ as
a quotient of a fiber bundle, all ``hands on''.
\par 
We start with the measure space $(\Omega,\calB(\Omega),\calI^k_\infty)$
where $1 \leq k \leq m-1$ are integers and $\calI^k_\infty$ is
the integral geometric measure described in \cite[2.10.5(1)]{GMT} and
\cite[5.14]{MATTILA}.
Note that it is not semifinite, \cite[3.3.20]{GMT}.
Thus, we replace it with its complete semifinite version
$(\Omega,\widetilde{\calB(\Omega)},\tilde{\calI}^k_\infty)$.
We let $\calE$ be the collection of $k$-dimensional submanifolds $M
\subset \Omega$ of class $C^1$ such that $\phi_M = \calH^k \hel M$ is
locally finite.
It follows from the Besicovitch Structure Theorem \cite[3.3.14]{GMT}
that $\calE$ is $\calN_{\tilde{\calI}^k_\infty}$-generating.
Now, for each $x \in \Omega$ we define $\calE_x = \calE \cap \{ M : x \in
M \}$ and we define on $\calE_x$ an equivalence relation as follows.
We declare that $M \sim_x M'$ if and only if
\begin{equation*}
  \lim_{r \to 0^+} \frac{\calH^k(M \cap M' \cap \bB(x,r))}{\balpha(k)r^k} 
= 1.
\end{equation*}
Letting $[M]_x$ denote the equivalence class of $M \in \calE_x$ (of which we may think as the germ of a $k$-rectifiable set at $x$), we
prove in \cite{BOU.DEP.21} that the underlying set of the cccc, lld, and decomposable version of the MSN
$(\Omega,\widetilde{\calB(\Omega)},\calN_{\tilde{\calI}^k_\infty})$ can be
taken to be
\begin{equation*}
  \hat{\Omega} = \{ (x,[M]_x ) : x \in \Omega \text{ and } M \in \calE_x \} \,.
\end{equation*}
This leads to an explicit description of the dual of
$\bL_1(\Omega,\widetilde{\calB(\Omega)},\tilde{\calI}^k_\infty)$ as
$\bL_\infty( \hat{X},\hat{\calA},\hat{\calN})$.
In particular case $k=n-1$ of interest for studying the dual of $\mathbf{SBV}(\Omega)$ we have:
\begin{Theorem*}
  The map $\Upsilon \colon \bL_\infty(\hat{\Omega}, \hat{\calA},
  \hat{\calI}^{n-1}_\infty) \to \bL_1(\Omega, \widetilde{\calB(\Omega)}, \tilde{\calI}_\infty^{n-1})^*$
  defined by
  \[
  \Upsilon(\bg)(\mathbf{f}) = \int_{\hat{\Omega}} g(f \circ p) \, d\hat{\calI}^{n-1}
  \]
  is an isometric isomorphism (where $p$ is as the last theorem of the previous section).
\end{Theorem*}
\par 
The method relies on a ``density 1'' theorem for $k$-rectifiable Borel sets and the fact that null sets in this context are purely $k$-unrectifiable.
The main point, of course, of this effort is that the quotient fiber bundle $\hat{\Omega}$ above is a ``localization'' in the sense that it is ``as small as possible'' in the role of a decomposable version of the integral geometric measure space that we started with.

\section{The dual of $\mathbf{SBV}$}
\label{sec.sbv}

With the quotient fiber bundle $\hat{\Omega}$ described in the previous section we can now state our ``optimal'' integral representation of members of the dual of $\mathbf{SBV}(\Omega)$, see \cite{DEP.BOU.22}.
The replacement of the Hausdorff measure by the integral geometric measure is due to the fact that we do not know an explicit description of the Radon-Nikod\'ymification of a Hausdorff measure (the problem is with the lack of a ``density 1'' theorem).
The reason why we may proceed with this replacement is that if $u$ is of special bounded variation then $\|Du\| << \calI^{n-1}_\infty$.

\begin{Theorem*}
  Let $\varphi \colon \mathbf{SBV}(\Omega) \to \R$ be a continuous linear
  functional. There are fields $f \in L_\infty(\Omega, \calB(\Omega),
  \calL^n)$, $g \in L_\infty(\hat{\Omega}, \hat{\calA},
  \hat{\calI}^{n-1}_\infty ; \R^n)$ and $h \in L_\infty(\Omega,
  \calB(\Omega), \calL^n ; \R^n)$ such that
  \begin{equation*}
    \varphi(u) = \int_\Omega fu \, d \calL^n + \int_{\hat{\Omega}} g \cdot (j_u \circ
    p) d \hat{\calI}^{n-1}_\infty + \int_\Omega h \cdot \nabla u \, d \calL^n
  \end{equation*}
  where the distributional gradient $Du$ of $u$ is decomposed into its
  Lebesgue part $\calL^n \hel \nabla u$ and its jump part $D^j u = j_u
  \, \calH^{n-1} \hel S_u$. Moreover, one can select $g$ in such a way
  that $g(x, [E]_x)$ is normal to $\rmTan^{n-1}([E]_x)$ for
  $\hat{\calI}^{n-1}_\infty$-almost all $(x, [E]_x) \in \hat{\Omega}$.
\end{Theorem*}

\section{Acknowledgements}

My thanks are due to {\sc David H. Fremlin}, not only for his inspiring treatise ``Measure Theory'' \cite{FREMLIN.I,FREMLIN.II,FREMLIN.III,FREMLIN.IV,FREMLIN.V.1,FREMLIN.V.2}, but also for many helpful conversations.
It is also my pleasure to record useful conversations with {\sc Francis Borceux} regarding category theory.
I am grateful to {\sc ZhiQiang Wang} for his careful reading of \cite{DEP.21} when it was still a manuscript and for his witty comments.
Last but not least, I am happily acknowledging my long and fruitful collaboration with {\sc Philippe Bouafia}.


\bibliographystyle{amsplain}
\bibliography{/Users/thierry/Documents/LaTeX/Bibliography/thdp}

\end{document}